\documentclass[twoside]{article}
\makeatletter
\def\ps@myheadings{
     \let\@oddfoot\@empty\let\@evenfoot\@empty
     \def\@evenhead{\thepage\hfil\MakeUppercase{\footnotesize\leftmark}\hfil}%
     \def\@oddhead{{\hfil\MakeUppercase{\footnotesize\rightmark}}\hfil\thepage}%
     \let\@mkboth\@gobbletwo
     \let\sectionmark\@gobble
     \let\subsectionmark\@gobble     }
\makeatother
\usepackage{amssymb,amsmath,amsfonts,amscd,amsthm}
\theoremstyle{plain}
\newtheorem{thm}{Theorem}[section]

\newtheorem{que}[thm]{Question}
\newtheorem{prob}[thm]{Problem}
\newtheorem{cor}[thm]{Corollary}
\newtheorem{prop}[thm]{Proposition}
\theoremstyle{definition}
\newtheorem{defn}[thm]{Definition}

\begin{document}
\pagestyle{myheadings} \markboth{Denise M. Halverson, Du\v{s}an
Repov\v{s}} {Survey on the Generalized R.~L.~Moore Problem}
\title{Survey on the Generalized R.~L.~Moore Problem}
\author{Denise M. Halverson and Du\v{s}an Repov\v{s}}
\date{\today}
\maketitle
\begin{abstract}
We give an updated extended survey of results related to the
celebrated unsolved generalized R.~L.~Moore problem.  In particular,
we address the problem of characterizing codimension one
manifold factors, i.e. spaces $X$ having the property that $X \times
\mathbb{R}$ is a topological manifold. A main part of the paper is
devoted to many efficient general position techniques, which have
been used to solve special cases of this problem.
\end{abstract}
\maketitle

\section{Introduction}

Let $\mathbb R^n$ denote the $n$-dimensional Euclidean space.
The Generalized R.~L.~Moore Problem asks: \medskip

\emph{Suppose that $G$ is a cellular (or
cell-like)
upper
semicontinuous
decomposition of $\mathbb{R}^n$, where $n\ge 3$.
Is then
$(\mathbb{R}^n/G) \times \mathbb{R}$
homeomorphic to $\mathbb{R}^{n+1}$?}
\medskip

\noindent This is a classical problem that has remained unsolved for
over sixty years.  One major importance of this problem is its
potential applications to manifold recognition problems such as the
famous Busemann Conjecture and Bing-Borsuk Conjecture --
see our recent survey \cite{HaRe}.

In 1980, Daverman published a most excellent survey on the
Generalized R.~L.~Moore Problem \cite{Daverman S}.  In the present
paper, we
extend this discussion to significant developments since that point
in time. Although we provide a few of the most relevant details from
the earlier period, we refer the reader to Daverman's survey for a
more thorough discussion. Our survey will focus on developments
after 1980, especially with respect to general position
strategies.

\section{Background}

The Generalized R.~L.~Moore Problem first emerged in the
investigation of manifold recognition type problems in the mid
1900's.  A decomposition $G$ of a Hausdorff space $S$ into compact
subsets is
\emph{upper semicontinuous} if and only if the
decomposition map $\pi: S \to S/G$ is closed. An early result of
R.~L.~Moore states \cite{Moore2}:

\begin{thm}[{\bf R.~L.~Moore Theorem}]
If $G$ is an upper semicontinuous decomposition of $\mathbb{R}^2$
into 
continua which do not
separate
$\mathbb{R}^2$, then the
decomposition space $\mathbb{R}^2/G$ is homeomorphic to
$\mathbb{R}^2$.
\end{thm}

Could this result be generalized to higher dimensions? Moore's
theorem was proved by appealing to a topological characterization of
the plane. Since no analogous characterizations for higher
dimensional manifolds existed at the time, generalizations of his
theorem to higher dimensions using the same approach were not
possible.

The first breakthrough in attacking higher dimensional problems came
during the 1950's when Bing developed a shrinkability criterion, which is used to define shrinkability in the following definition:

\begin{defn}
An upper semicontinuous decomposition $G$ of a space $X$ is said to
be \emph{shrinkable} if and only if $G$ satisfies the following shrinkability criterion: for each open cover $\mathcal{U}$ of $X/G$
and for each open cover $\mathcal{V}$ of $X$ there exists a
homeomorphism $h$ of $X$ onto itself satisfying:
\begin{enumerate}
\item  To each $x \in X$ there corresponds $U \in \mathcal{U}$ such
that $\{ x, h(x) \} \subset \pi^{-1}(U)$; and
\item  To each $g \in G$ there corresponds $V \in \mathcal{V}$ such
that $h(g) \subset V$.
\end{enumerate}
\end{defn}

\noindent With various minor hypotheses placed on the source space,
it can be determined whether or not the decomposition map $\pi: X
\to X/G$ is a near-homeomorphism, i.e. $\pi$ can be approximated by homeomorphisms.  One of the most general cases is addressed in the following theorem (see
\cite{Edwards 2,EG,MV}):

\begin{thm}
Suppose $G$ is a usc decomposition of a complete metric space $X$.
Then the decomposition map $\pi: X \to X/G$ is a near-homeomorphism
if and only if $G$ is shrinkable.
\end{thm}

\noindent
With this and similar theorems,
Bing's  shrinkability criterion
opened the door wide for exploring many examples of
decomposition spaces.

Early on, the cellularity property was investigated for its
potential as a condition on decomposition elements which might imply
a similar result as Moore's Theorem for decompositions of
$\mathbb{R}^3$. A subset $X$ of an $n$-manifold $M$ is said to be
\emph{cellular} in $M$ if $M$ contains a family of $n$-cells $\{C_i
\ | \ i = 1,2, \ldots \}$ such that $C_{i+1} \subset \text{int }
C_i$ and $X = \bigcap C_i$. In particular, could it be true that if
$G$ is an upper semicontinuous decomposition of $\mathbb{R}^3$ into
cellular sets, then is the decomposition space $\mathbb{R}^3/G$ is
homeomorphic to $\mathbb{R}^3$?  In 1957, Bing \cite{bing1}
constructed the
Dogbone space which was realized by an upper semicontinuous cellular
decomposition of $\mathbb{R}^3$ but failed to be homeomorphic to
$\mathbb{R}^3$, thereby demonstrating that the answer
to the previous question is no. Many examples of decompositions of
manifolds whose elements are cellular, or even more generally
cell-like, that do not generate manifolds, have now been discovered
(cf.  \cite{Daverman book} for a catalog of examples).

Shortly after the construction of his Dogbone space, Bing \cite{bing2}
discovered
a very surprising result: the product of the Dogbone space with the
real line is homeomorphic to $\mathbb{R}^4$. The
inherent "tangling" of the decomposition elements preventing
shrinking becomes sufficiently ``unraveled'' upon taking the product
with $\mathbb{R}$ so that the decomposition elements can be
simultaneously shrunk and the shrinkability criterion satisfied. Is this
always the case? Or could elements be so tangled, that even within a
product of $\mathbb{R}$ there is insufficient room to obtain the
desired shrinking of elements. Thus emerged the {\it Generalized R.~L.~Moore Problem}:

\begin{prob}[{\bf Generalized R.~L.~Moore Problem}]
If $G$ is a cellular upper
semicontinuous decomposition of
$\mathbb{R}^n$, is $\mathbb{R}^n/G \times \mathbb{R}$ homeomorphic
to $\mathbb{R}^{n+1}$?
\end{prob}

Later, it was determined that the condition ``cellular'' may be more
appropriately replaced with ``cell-like''. An \emph{absolute
neighborhood retract (ANR)} is a locally contractible Peano
continuum. A compact subset $X$ of an ANR $Y$ is \emph{cell-like} if
$X$ is contractible in every neighborhood of itself.  Cellular sets
are necessarily cell-like, but the converse is not true. Moreover, a
set being cell-like is inherent to the set itself, whereas
cellularity depends how the set is embedded in the ambient space. It
is the cell-like condition on an upper semicontinuous decomposition
of a manifold that is sufficient to imply that the decomposition
space is a homology manifold \cite{Repovs1}--\cite{Repovs3}.

A few examples of some early results concerning the Generalized
R.~L.~Moore Problem include the following: \medskip

\noindent The product $(E^n/G) \times E^1$ is homeomorphic to
$E^{n+1}$ in the following cases:
\begin{itemize}
\item $G$ consists of a single nondegenerate element that is an arc
(Andrews-Curtis \cite{Andrews-Curtis}).
\item $G$ consists of a single nondegenerate element that is a cell
(Bryant \cite{Bryant 1, Bryant 2}).
\item $G$ is an usc decomposition into points and a countable
collection of arcs (Gillman-Martin \cite{Gillman-Martin}).
\item $G$ is an usc decomposition into points and a null sequence of
cells (Meyer \cite{Meyer}).
\end{itemize}

\noindent These examples provide a sense of the types of problems
that were most accessible through shrinking theorems.

\section{Characterizations of Manifolds}

Essential to making progress on the Generalized R.~L.~Moore Problem
are effective characterizations of manifolds.  In dimensions $1$ and
$2$, characterizations of manifolds are relatively simple. For
example, a space is homeomorphic to $S^1$ if and only if it is a
nondegenerate locally connected continuum (compact connected
Hausdorff space) that is separated by no single point, but is
separated by any pair of points \cite{Moore1, Wilder}. The {\it Kline Sphere
Characterization Theorem} states that a space is homeomorphic to
$S^2$ if and only if it is a nondegenerate locally connected metric
continuum that is separated by any simple closed curve but not
separated by any pair of points \cite{bing}.

When the Generalized Moore Problem was first posed,
characterizations of manifolds useful in addressing this problem in
higher dimensions were lacking. However, more breakthroughs were
initiated in the 1970's when Cannon proved the double suspension
problem and Edwards discovered the proof of the famous Cell-like
Approximation Theorem \cite{Daverman book, Daverman-Halverson,
Edwards 2} which characterizes manifolds of dimension $n
\geq 5$.  Before stating the Cell-like Approximation Theorem, we
shall review some of the essential properties of manifolds.

\subsection{Manifolds are generalized manifolds}

Throughout this paper we shall work with singular homology
with integer coefficients.
A {\it homology $n$-manifold} is a locally compact Hausdorff space
$X$ such that for every $x \in X$, $H_i(X, X-\{x\}) \cong
H_i(\mathbb{R}^n, \mathbb{R}^n-\{0\})$ for all integers $i \geq 0$.
A {\it homology $n$-manifold with boundary} is a locally compact
Hausdorff space $X$ such that either $H_i(X, X-\{x\}) \cong
H_i(\mathbb{R}^n, \mathbb{R}^n-\{0\})$ or $H_i(X, X-\{x\}) \cong 0$
for all integers $i \geq 0$.   The points $x$ such that $H_i(X,
X-\{x\}) \cong 0$ for all integers $i \geq 0$ are called the
boundary points (cf. \cite{Bredon,Mitchell}).

A {\it Euclidean neighborhood retract} (ENR) is a
space $X$ that embeds in $\mathbb{R}^n$ so that $X$ is a retract of
a neighborhood of itself in $\mathbb{R}^n$. ENR's are the finite-dimensional
ANR's. A \emph{generalized $n$-manifold} is an ENR
homology $n$-manifold. Clearly manifolds are generalized manifolds.
However, amongst many other examples, Bing's Dogbone space
demonstrates that not all resolvable generalized manifolds are
manifolds.

\subsection{Manifolds are resolvable}

A map $f:Y \to X$ is \emph{cell-like} if for every $x \in X$,
$f^{-1}(x)$ is cell-like
(for more on this important class of maps see the survey \cite{MiRe}).
A space $X$ is said to be
\emph{resolvable} if there is a manifold $M$ and a surjective map
$f: M \to X$ which is cell-like. Such a map $f$ is said to be a
\emph{resolution} of $X$.  In this case, $X$ is said to be
\emph{resolvable}.  Clearly, manifolds are resolvable.

A large class of generalized 3-manifolds is known to be resolvable
\cite{Thickstun 1,Thickstun 2}.
For admissible generalized manifolds of dimension $n\geq 4$, the
question of whether or not the space is resolvable can be determined
by a locally defined number called the {\it Quinn index} (cf.  \cite{Quinn
2}).

\begin{defn}
Let $Y$ be a homology $n$-manifold.
\begin{itemize}
\item If $\partial Y = \emptyset$, then $Y$ is said to be {\it admissible
to resolution theory} provided $Y$ is a locally compact, connected,
finite-dimensional, separable, metrizable ANR.
\item If $\partial Y \ne \emptyset$, the $Y$ is said to be
{\it admissible to resolution theory} provided both $Y-\partial Y$ and
$\partial Y$ are admissible homology manifolds with empty boundary.
\end{itemize}
\end{defn}

\begin{thm} [{\bf Local Index Theorem}]
Let $n\geq 4$ and let $Y$ be an admissible homology $n$-manifold.
There is a local index $i(Y) \in (1+8\mathbb{Z})$ defined, which has
the following properties:
\begin{enumerate}
\item For every non-empty, open $U \subset Y$, $i(U) = i(Y)$;
\item If $\partial Y \ne \emptyset$, then $i(\partial Y) = i(Y)$;
\item If $X$ is admissible, then $i(X \times Y) = i(X) \times i(Y)$; and
\item If dim $Y \geq 5$, or if dim $Y = 4$ and $\partial Y$ is
either empty or a $3$-manifold, then $Y$ is resolvable if and only
if $i(Y)=1$.
\end{enumerate}
\end{thm}

It is unknown if all generalized $n$-manifolds are
resolvable for $3\leq n \leq 6$ (see \cite{Repovs2}).
Bryant, Ferry, Mio, and
Weinberger \cite{BFMW}
demonstrated the existence of nonresolvable generalized
manifolds in dimensions $n \geq 7$.

\subsection{Manifolds have general position properties}

General position properties deal with the ability to separate mapped
in objects by small adjustments, based on their dimension.  We shall
see that these types of properties are very useful in detecting both
manifolds and codimension one manifold factors.

The disjoint disks properties are the most basic of all types of
general position properties.  Let $D^j$ denote a $j$-cell.  A space
$X$ is said to satisfy the \emph{$(m,k)$-disjoint disks property}
($(m,k)$-DDP) if any two maps $f:D^m \to X$ and $g:D^k \to X$ can be
approximated by maps with disjoint images. Certain $(m,k)$-disjoint
disks properties have special names:

\begin{itemize}
\item The $(0,2)$-DDP is called the \emph{disjoint point-disk property}
(DPDP).
\item The $(1,1)$-DDP is called the \emph{disjoint arcs property} (DAP).
\item The $(1,2)$-DDP is called the \emph{disjoint arc-disk property}
(DADP).
\item The $(2,2)$-DDP is called the \emph{disjoint disks property}
(DDP).
\end{itemize}

Inherent to manifolds are the general position properties with
respect to the appropriate dimensions:

\begin{prop}
An $n$-manifold has the $(m,k)$-DDP if $m+k+1 \leq n$.
\end{prop}

\noindent However, general position for an arbitrary resolvable
generalized manifolds may be far more restrictive.  The following
important result for can be found in \cite[Proposition
26.3]{Daverman book}:

\begin{prop}
A generalized $n$-manifold, for $n \geq 3$, has the $(1,1)$-DDP.
\end{prop}

\noindent Beyond this, little is guaranteed.  Daverman and Walsh
\cite{Daverman-Walsh}
poignantly demonstrated this fact in their construction of ghastly
spaces, resolvable generalized manifolds of dimension $n \geq 3$
that fail to have the $(0,2)$-DDP.

\subsection{The Cell-like Approximation Theorem}

In the 1970's, major breakthroughs occurred in the manifold
recognition problem. In proving the famous double suspension
problem, Cannon first recognized the significance of the disjoint
disks property in characterizing high dimension manifolds
\cite{Cannon, Cannon 1}. His insight was affirmed by the monumental
proof of the Edwards Cell-like Approximation Theorem (cf.
\cite{Daverman book} for the case $n>5$ and
\cite{Daverman-Halverson} for the case $n=5$):

\begin{thm}[{\bf Cell-like Approximation Theorem}]
For all $n\geq 5$, topological $n$-manifolds are precisely
the resolvable generalized $n$-manifolds that have the disjoint
disks property.
\end{thm}

It is well known that not all resolvable generalized manifolds of
dimension $n \geq 5$ have the DDP (cf. \cite{Daverman book}). Thus,
as indicated previously, not all resolvable generalized manifolds
are manifolds. Daverman and Repov\v s introduced
an appropriate analogue of
DDP in dimension $3$
which yields analogous results \cite{Repovs1}-\cite{Repovs3} (cf. also \cite{DaTh}),
whereas in dimension $4$ we still do not have a good replacement for
DDP (there exist so far  only some taming theorems -- cf. \cite{BDVW}).

\section{Manifold Factors}

Convenient for our discussion is the terminology of manifold
factors.  A {\it codimension $k$ manifold factor} is a space $X$
such that $X \times \mathbb{R}^k$ is a manifold. Thus a {\it
codimension one manifold factor} is a space $X$ such that $X \times
\mathbb{R}$ is a manifold.

Daverman has demonstrated the following intriguing result in
\cite{Daverman 1}:

\begin{thm}
If $X$ is a resolvable generalized manifold of finite dimension $n
\geq 3$, then $X \times \mathbb{R}^2$ is a manifold.
\end{thm}

\noindent It follows that all finite-dimensional resolvable
generalized manifolds are codimension $k$ manifold
factors for $k
\geq 2$.  The question of whether or not all finite-dimensional
resolvable generalized manifolds are codimension one manifold
factors remains unsolved.  In the terminology of manifold factors, a
more general statement of the Generalized R.~L.~Moore Problem is the
{\it Product with a Line Problem}:

\begin{prob}[{\bf Product with a Line Problem}]
Characterize $n$-dimensional spaces $X$ such that $X \times
\mathbb{R}$ is an $(n+1)$-dimensional manifold.
\end{prob}

\noindent A necessary condition for a space to be a manifold factor
is that the space is a resolvable generalized manifold \cite{Quinn
2}. In the case of spaces of dimension $n=3$, the verification that
a space is a codimension one manifold factor is still generally done
by shrinking arguments. However, more flexible methods of
verification in the case of $n \geq 4$ are possible due to the
characterization of $n$-manifolds, $n \geq 5$, provided by Edward's
cell-like approximation theorem.

\section{Detecting Codimension One Manifold Factors with General Position Properties}

There are several general position properties that are useful in
detecting codimension one manifold factors.  Each one is designed so
that if $X$ has the said property, then $X \times \mathbb{R}$ has
the disjoint disks property. The three main general position
properties that we shall discuss, in order of strength in detecting
codimension one manifold factors are: the disjoint arc-disk
property, the disjoint homotopies property, and the disjoint
concordances property or equivalently the disjoint topographical
maps property. In each case, conditions known to implicate the said
property are also included.


\subsection{Disjoint Arc-Disk Property}

Perhaps the simplest analog to the disjoint disks property in one
dimension lower is the disjoint arc-disk property DADP, mentioned
earlier as the (1,2)-DDP.  Daverman has shown the following
\cite{Daverman 1}:

\begin{thm}
If $X$ is a resolvable generalized manifold with the DADP, then $X$
is a codimension one manifold factor.
\end{thm}

\noindent Resolvable generalized manifolds that possess the disjoint
arc-disk property include spaces that arise from $(n-3)$-dimensional
or closed $(n-2)$-dimensional decompositions \cite{CaBrLa}. Examples
include decompositions of $n$-manifolds arising from a classical
defining sequence \cite[Proposition 9.1]{Daverman book} and
decompositions of $\mathbb{R}^{n \geq 4}$ into convex sets
\cite{HaRe3}.

However, not all manifolds known to be codimension one manifold
factors have the disjoint arc-disk property.  Examples include the
totally wild flow  \cite{Cannon-Daverman2}, the Daverman-Walsh
ghastly spaces \cite{Daverman-Walsh} and the
$k$-ghastly spaces \cite{Halverson 1}. Thus a more effective general
position property is desired.

\subsection{Disjoint Homotopies Property}

The disjoint homotopies property has proven to be an effective
detector of codimension one manifold factors in almost every known
case.

Let $D=I=[0,1]$.  A space $X$ has the \emph{disjoint homotopies
property (DHP)} if every pair of path homotopies $f,g: D \times I
\to X$ can be approximated by disjoint homotopies, i.e.,
approximating maps $f',g': D \times I \to X$ so that $f_t(D) \cap
g_t(D) = \emptyset$ for all $t \in I$.  The following theorems are
proven in \cite{Halverson 1}:

\begin{thm}  If $X$ is a locally compact ANR with DHP, then $X
\times \mathbb{R}$ has DDP. \end{thm}

\begin{cor}
If $X$ is a resolvable generalized $n$-manifold, where $n\geq 4$,
having DHP, then $X$ is a codimension one manifold factor.
\end{cor}

Strategies for obtaining approximating disjoint homotopies include:
\begin{enumerate}
\item Reimaging, or adjusting the image set.
\item Realigning the levels within the domain.
\item Reparameterizing the levels.
\end{enumerate}

\noindent Moreover, it is sufficient to prove that two path
homotopies can be approximated by disjoint homotopies in the case
that one of the path homotopies is constant.

\begin{thm} \label{jiggle}
If $X$ is an ANR with the property that any constant path homotopy
together with any arbitrary path homotopy can be approximated by
disjoint homotopies, then $X$ has DHP.
\end{thm}

\noindent These are the strategies employed to demonstrate that DHP
is implied for spaces possessing one of the the next three related
properties.

\subsubsection{Plentiful 2-manifolds property}

A space $X$ has the \emph{plentiful $2$-manifolds property} (P2MP)
if each path $\alpha: I \to X$ can be approximated by a path
$\alpha':I \to N\subset X$ where $N$ is a $2$-manifold embedded in
$X$.

\begin{thm}
Suppose $X$ is a generalized $n$-manifold, $n\geq 4$, $g:D \times I
\to X$ and$f:D \times I \to N \subset X$ where $N$ is a $2$-manifold
embedded in $X$.  Then $f$ and $g$ can be approximated by disjoint
homotopies.
\end{thm}

\begin{cor}
If $X$ is a resolvable generalized $n$-manifold, $n \geq 4$, with
P2PM, then $X$ is a codimension one manifold factor.
\end{cor}

Examples of spaces that have the plentiful $2$-manifolds property
are decomposition spaces resulting from a nested defining sequences
of thickened $(n-2)$-manifolds, spaces that arise from closed
$(n-2)$-dimensional decompositions, and certain $k$-ghastly spaces
for $2 < k < n$ (cf. \cite{Halverson 1}).

\subsubsection{0-Stitched disks property}

The maps of $f,g:D^2 \to X$ are said to be \emph{$0$-stitched}
provided that there are $0$-dimensional $F_\sigma$ sets $A$ and $B$
contained in the interior of $D^2$ such that $f(D^2-A) \cap g(D^2-B)
= \emptyset$.  We say that $f$ and $g$ are \emph{$0$-stitched along
$A$ and $B$}.  If $Y$ and $Z$ are sets in $D^2$ missing $A$ and $B$
respectively, then we say that $f$ and $g$ are \emph{$0$-stitched
away from $Y$ and $Z$}.

A space $X$ has the \emph{$0$-stitched disks property} if any two
maps $f,g:D^2 \to X$ can be approximated  by maps $f',g':D^2 \to X$
such that $f'$ and $g'$ are $0$-stitched along sets $0$-dimensional
$F_\sigma$-sets $A$ and $B$ and away from infinite $1$-skeleta
$(K_j^\infty)^{(1)}$, $j=1,2$, of $D^2$ such that
$f'|_{(K_1^\infty)^{(1)}} \cup g'|_{(K_2^\infty)^{(1)}}$ is $1-1$.

\begin{thm}
If $X$ has the $0$-stitched disks property, then $X$ has DHP.
\end{thm}

There are many examples of resolvable generalized manifolds $X =
M/G$ in which the $0$-stitched disks property can be easily
verified. Suppose that $M$ is a manifold with a sequence of
triangulations $\{K_i\}$ such that mesh$(K_i) \to 0$ and  $G$ is a
usc decomposition arising so that the nondegeneracy set that misses
the infinite $1$-skeleta $(K_j^\infty)^{(1)}$ and meets the infinite
$2$-skeleta $(K_j^\infty)^{(2)}$ in a $0$-dimensional
$F_\sigma$-set. Then the resulting decomposition space will have the
$0$-stitched disks property. These conditions are generally easily
imposed in many constructions of resolvable generalized manifolds
that arise from defining sequences.

\subsubsection{The Method of Delta-Fractured Maps}

A map $f:D \times I \to X$ is said to be \emph{$\delta$-fractured}
over a map $g:D \times I \to X$ if there are pairwise disjoint balls
$B_1, B_2, \ldots, B_m$ in $D \times I$ such that:
    \begin{enumerate}
        \item $diam(B_i) < \delta$;
        \item $f^{-1}(im(g)) \subset \bigcup_{i=1}^m int(B_i)$; and
        \item $diam(g^{-1}(f(B_i))) < \delta$.
    \end{enumerate}

\begin{thm}
If $X$ is an ANR that has the property  for any pair of path
homotopies $f,g: D\times I \to X$, where $g$ is a constant path
homotopy, and $\delta
> 0$ their are approximations $f', g': D \times I \to X$ of $f$ and $g$, respectively, such that
$f'$ is $\delta$-fractured over $g'$, then $X$ has DHP.
\end{thm}

The strength of the method of $\delta$-fractured maps is manifest in
its application  to prove that certain $2$-ghastly spaces have DHP
(cf. \cite{Halverson 2}).

\subsection{Disjoint Concordance Property and the Disjoint Topographies Property}

Although the disjoint homotopies properties has proven extremely
useful in detecting codimension one manifold factors, it is still
unknown whether it is a necessary condition on codimension one
manifold factors of dimension $n \geq 4$. However,  the disjoint
concordance property was shown in \cite{Daverman-Halverson}  to be
both a necessary and sufficient condition on resolvable generalized
manifolds of dimension $n \geq 4$ that are codimension one manifold
factors.

A \emph{path concordance} in a space $X$ is a map $F:D \times I \to
X \times I$ (where $D=I=[0,1]$) such that $F(D \times e) \subset X
\times e, e \in \{0,1\}.$  A metric space ($X,\rho$) satisfies the
\emph{disjoint path concordances property (DCP)} if, for any two
path homotopies $f_i:D \times I \to X$ ($i=1,2$) and any
$\varepsilon > 0$, there exist path concordances $F'_i: D \times I
\to X \times I$ such that
\begin{center}
$F'_1(D \times I) \cap F'_2(D \times I) = \emptyset$
\end{center}
and $\rho (f_i, \text{proj}_X F'_i) < \varepsilon$.

\begin{thm}[Daverman and Halverson \cite{Daverman-Halverson 2}] Suppose $X$ is a locally
compact, metric ANR with DAP. Then $X$ has DCP if and only if $X
\times \mathbb{R}$ has DDP.
\end{thm}

\begin{cor}
A resolvable generalized manifold $X$ of dimension $n\geq 4$ is a
codimension one manifold factor if and only if $X$ has DCP.
\end{cor}

The problem with the disjoint concordances property is that in its
raw form, it has not had a great deal of utility.  However, the
disjoint topographies property, a condition equivalent to DCP but
 having more of the flavor of the disjoint homotopies
property,  has much more potential.  The advantage of the disjoint
topographies property over DHP is that not only does it allow for
the same strategies of reimaging, realigning, and reparameterizing,
but it also allows for the change in the shape of the levels.

A \emph{topography $\Upsilon$ on $Z$} is a partition of $Z$ induced
by a map $\tau: Z \to I$.  The \emph{$t$-level of $\Upsilon$} is
given by $$ \Upsilon_t = \tau^{-1}(t).$$

A \emph{topographical map pair} is an ordered pair of maps $(f,
\tau)$ such that $f:Z \to X$ and $\tau: Z \to I$.  The map $f$ will
be referred to as the \emph{spatial map} and the map $\tau$ will be
referred to as the \emph{level map}. The topography associated with
$(f, \tau)$ is $\Upsilon$, where $ \Upsilon_t = \tau^{-1}(t).$

Note that a homotopy $f: Z\times I \to X$ has a naturally associated
topography, where $\tau: Z\times I \to I$ is defined by $\tau(x,t) =
t$. In particular, we may view $f: Z\times I \to X$ as being
equivalent to $(f,\tau)$ and we shall refer to $(f,\tau)$ as the
\emph{natural topographical map pair associated with $f$}.

Suppose that for $i=1,2$, $\Upsilon^i$ is a topography on $Z_i$
induced by $\tau_i$ and $f_i: Z_i \to X$.  Then $(f_1, \tau_1)$ and
$(f_2, \tau_2)$ are \emph{disjoint topographical map pairs} provided
that for all $t \in I$,
$$f_1( \Upsilon^1_t) \cap f_2(\Upsilon^2_t) = \emptyset.$$
A space $X$ has the \emph{disjoint topographies property (DTP)} if
any two topographical map pairs $(f_i, \tau_i)$ ($i=1,2$), where
$f_i: D^2 \to X$, can be approximated by disjoint topographical map
pairs.

\begin{thm} [Halverson and Repov\v{s} \cite{HaRe2}]
An ANR $X$ has the disjoint topographies property if and only if $X
\times \mathbb{R}$ has DTP.
\end{thm}

\begin{cor}
A resolvable generalized manifold $X$ of dimension $n\geq 4$ is a
codimension one manifold factor if and only if $X$ has DTP.
\end{cor}

The following ribbons properties have analogs to the special
properties defined for DHP.  The crinkled ribbons properties are a
generalization of the plentiful $2$-manifolds property.  The fuzzy
ribbons property is a generalization of the method of
$\delta$-fractured maps (cf. \cite{HaRe2}).

\subsubsection{The Crinkled Ribbons Properties}

A generalized $n$-manifold $X$ has the \emph{crinkled ribbons
property (CRP)} provided that any constant homotopy $f: K \times I
\to X$, where $K$ is a $1$-complex can be approximated by a map
$f':K \times I \to X$ so that:
\begin{enumerate}
\item $f'(K \times \{0\}) \cap f'(K \times \{1\}) = \emptyset$; and
\item dim$(f'(K \times I))\leq n-2$.
\end{enumerate}

\begin{thm}
If $X$ is a resolvable generalized $n$-manifold, $n \geq 4$, with
the crinkled ribbons property, then $X$ is a codimension one
manifold factor.
\end{thm}

A generalized $n$-manifold $X$ has the \emph{twisted crinkled
ribbons property (CRP-T)} provided that any constant homotopy $f: D
\times I$ can be approximated by a map $f':D \times I$ so that:
\begin{enumerate}
\item $f'(D \times \{0\} ) \cap f'(D \times \{1 \} )$ is a finite set of
points; and
\item dim$(f'(D \times I))\leq n-2$.
\end{enumerate}


\begin{thm}
If $X$ is a generalized $n$-manifold of dimension $n \geq 4$ having
 the twisted crinkled ribbons property and the property that points
are $1$-LCC embedded in $X$, then $X$ is a codimension one manifold
factor.
\end{thm}

One should note that not all generalized manifolds of dimension $n
\geq 4$ have the property that points are $1$-LCC embedded.  For
example, the Daverman-Walsh $2$-ghastly spaces are resolvable
generalized manifolds that do not have the $(0,2)$-DDP, and hence
cannot satisfy the condition that points are $1$-LCC embedded
\cite{Daverman-Walsh}.

One application of the CRP is the result that if $X$ is a resolvable
generalized locally spherical $n$-manifold, $n\geq 4$, then $X$ is a
codimension one manifold factor \cite{HaRe2}.  Although this result
was initially shown using shrinking arguments \cite{Daverman 2,
Daverman book}, it had not been proven previously using general
position techniques.

\subsubsection{The Fuzzy Ribbons Property}

Because of the freedom in restructuring the levels of the
topographies to obtain DTP conditions, the $\delta$-control in the
method of $\delta$-fractured maps is not required.  The analogous
definition of $\delta$-fractured maps in the setting of
topographical map pairs is as follows:

Let $K$ be a $1$-complex.  A topographical map pair $(f, \tau)$ is
in the \emph{$\mathcal{K}$ category} if $f:K \times I \to X$ and
$\tau: K \times I \to I$ so that $K \times \{e\} \subset
\tau^{-1}(e)$ for $e=0,1$. We denote $(f, \tau) \in \mathcal{K}$. A
topographical map pair $(f,\tau)$ is in the \emph{$\mathcal{K}_c$
category} if
\begin{enumerate}
\item $(f,\tau) \in \mathcal{K}$;
\item $f:K \times I \to X$ is a constant
homotopy; and
\item $(f,\tau)$ is the natural topographical map pair associated with $f$.
\end{enumerate}

Let $(f_i, \tau_i) \in \mathcal{K}$ be such that $f_i: K_i \times I
\to X$ and $\tau_i: K_i \times I \to I$.  Then $(f_2, \tau_2)$ is
said to be {\it fractured} over a topographical map pair $(f_1,
\tau_1)$ if there are disjoint balls $B_1, B_2, \ldots, B_m$ in $K_2
\times I$ such that:
    \begin{enumerate}
        \item $f_2^{-1}(im(f_1)) \subset \bigcup_{j=1}^m int(B_i)$; and
        \item $\tau_1 \circ f_1^{-1}\circ f_2(B_i) \ne I$.
    \end{enumerate}

A space $X$ has the \emph{fuzzy ribbons property (FRP)} provided
that for any topographical map pairs, $(f_1, \tau_1) \in
\mathcal{K}_c$ and $(f_2, \tau_2) \in \mathcal{K}$, and $\varepsilon
> 0$ there are maps $\tau'_i$ and $\varepsilon$-approximations $f'_i$ of  $f_i$  so that
$(f'_2, \tau'_2)$ is fractured over $(f'_1, \tau'_1)$.

\begin{thm}
If $X$ is a generalized $n$-manifold of dimension $n \geq 4$ having
the fuzzy ribbons property, then $X$ is a codimension one manifold
factor.
\end{thm}

Certain $2$-ghastly spaces satisfy the FRP, such as those discussed
in \cite{Halverson 2}. The same type of arguments apply, however
less attention to control is needed to satisfy the FRP.


\section{Epilogue}

We list a few interesting unsolved problems:

\begin{que}
   If $G$ is an $(n-2)$-dimensional cell-like decomposition of an
$n$-manifold $M$, where $n \geq 4$, is $M/G$ a codimension one
manifold factor?
\end{que}

\begin{que}
Is every finite-dimensional resolvable generalized manifold of
dimension $n\geq 4$ a codimension one manifold factor?
\end{que}

Beginning in 1942, Busemann \cite{Busemann1, Busemann2} developed
the notion of a $G$-space as a way of putting a Riemannian like
geometry on a metric space (and also in an attempt to obtain a
"synthetic description" of Finsler's spaces \cite{Finsler}). A {\it
Busemann $G$-space} is a metric space that satisfies four basic
axioms on a metric space. These axioms infer that Busemann
$G$-spaces are homogeneous geodesic spaces with the property that
small metric balls have a cone structure.

Busemann  \cite{Busemann3}
conjectured
that
every $n$-dimensional Busemann $G$-space ($n \in \mathbb{N}$) is a
topological $n$-manifold. This conjecture has been proven true
for dimensions $n \leq 4$ \cite{Busemann2, Krakus, Thurston}.
The Busemann Conjecture is also known to be true in all dimensions under the additional hypothesis that the Aleksandrov curvature is bounded either from below or from above \cite{Berestovskii 1, Berestovskii 2}.

In the general setting, the solution to the Busemann Conjecture is determined by the answer to the following question:

\begin{que} \label{Q}
Are small metric spheres in $n$-dimensional Busemann
$G$-spaces ($n \in
\mathbb{N}$) codimension one manifold factors?
\end{que}

$G$-homogeneous Busemann $G$-spaces are spaces in which the 
cone structure of the small metric balls are stable near their cone point. 
As a possible clue to the answer to Question \ref{Q}, it has been shown
that in the case of $G$-homogeneous Busemann $G$-spaces, small metric 
spheres are homogeneous. Moreover, these spaces need not 
have Aleksandrov curvature bounded either from below or from above \cite{BeHaRe}. 
It is unknown whether all Busemann $G$-spaces are $G$-homogeneous.

As for the prognosis, we believe that the first problem to tackle
should possibly be Question 6.1 since it appears the most tractable.
On the other hand, as our paper shows,  the
Generalized R.~L.~Moore Problem remains a
formidable question  - in spite of the great amount of work done in
the last half of the century - and it will probably occupy
generations to come.

\section*{Acknowledgements}
The authors were supported in part by the Slovenian Research Agency grants
BI-US/09-12/004,
P1-0292-0101,
J1-2057-0101
and
J1-4144-0101.
This paper was presented at the Conference on Computational and Geometric Topology
(Bertinoro, June 17-19, 2010).
We thank the organizers for the invitation and hospitality.

\end{document}